\documentstyle[12pt]{article} 

\begin{document}

\newcommand{\e}{\epsilon}
\renewcommand{\a}{\alpha}
\renewcommand{\b}{\beta}
\def\k{\kappa}
\newcommand{\p}{\partial}
\renewcommand{\phi}{\varphi}

\newtheorem{theorem}{Theorem}

\title{Nonlocal Dynamics of Passive Tracer Dispersion with Random 
Stopping}

\author{Hongjun Gao$^1$, James R. Brannan$^2$ and Jinqiao Duan$^2$ 
\\ 
\\1. Laboratory of Computational Physics\\ 
Institute of Applied Physics and Computational Mathematics\\ 
Beijing, 100088, China. 
\\ 
\\2. Department of Mathematical Sciences\\  
Clemson University, Clemson, South Carolina 29634, USA. \\ 
E-mail: duan@math.clemson.edu   }
 
      
\maketitle

\begin{abstract}

 We   investigate  the nonlocal behavior
of passive tracer dispersion with random stopping at various sites
in fluids.
This kind of dispersion processes is modeled by an integral
partial differential equation, i.e., an advection-diffusion equation
with a memory term.  We have shown the exponential decay of the
passive tracer concentration, under suitable conditions for the
velocity field and the probability distribution of random stopping time.

\bigskip

                 
\end{abstract}

\section{Introduction}

The discharge of pollutants  into coastal
seas or rivers is common but harmful.
The  pollutants may be held up or trapped 
at some sites
during the process of their dispersion.
For the benefit  of  
better environment, it is important to understand the
dynamics of such  passive tracers.

The Eulerian approach for studying passive tracer dispersion attempts to
understand the evolution of tracer concentration profile as a continuous 
field
quantity  \cite{Clark} \cite{Schnoor},  which ultimately, at large 
times,
satisfies the advection-diffusion equation.   

When there are side branches, traps or some special sites
in  a shear flow, where passive tracers
may be held-up or arrested and the stopping times of tracers at these 
traps
are random, 
the tracer concentration profile $C(x,t)$ then satisfies a
nonlocal transport equation derived by Young \cite{Young2} 
\begin{eqnarray}
  C_t + U(x) C_x     
  =  D C_{xx} -  \p_t \int_0^t  k (t-s) C(x,s) ds     
  \label{eqn}
\end{eqnarray} 
where  $D > 0$ is the diffusion constant and $U(x)$ is
the smooth fluid velocity field. The kernel $k(t) > 0$ is proportional 
to the probability distribution of tracer stopping times
at traps. 
  
Young \cite{Young2} has shown that the anomalous diffusion 
phenomenon can occur under certain conditions on velocity field
$U(x)$ and kernel $k(t)$. 
 
In this paper, we study the long-time behavior of the concentration 
profile.  We have shown the exponential decay of the
passive tracer concentration, under suitable conditions for the
velocity field and the probability distribution of random stopping time.

\section{Exponential Decay with \\
Increasing Random Stopping Distribution}
 
 We rewrite equation (\ref{eqn}) as
\begin{eqnarray}
  C_t + U(x) C_x     
  =  D C_{xx} - k(0) C - \int_0^t  k'(t-s) C(x, s) ds ,    
  \label{eqn2}
\end{eqnarray} 
where $x$ varies on the real line  $R = (-\infty, \infty)$.
We   impose an appropriate initial condition
\begin{eqnarray}
 C(x, 0) = C_0(x ).
 \label{IC}
\end{eqnarray}
 
We  denote $L^2 (R)$ as   the standard
  space of square-integrable  
functions. The usual mean-square norm in this
space is denoted as $\| \cdot \|$. 
Moreover, $H^1(R)$ denotes the usual Sobolev space.
Note that from \cite{Adams}, p. 57, 
$H^1(R) = H^1_0(R)$. We also denote 
$L^1(R^+)$ as the space of integrable functions
on positive real line $R^+ = (0, \infty)$.
Moreover, $C^1 (R^+)$ is the space of 
continuously differentiable (smooth) functions on $R^+$.

 In order to study the evolution of the
 concentration $C(x,t)$, we should  estimate the
 nonlocal convolution integral in (\ref{eqn2}), and thus 
  we need the following   result  \cite{sta}.

\bigskip

\noindent
{\bf Lemma.}  
{\em If $k\in L^1(R^+)$ is a positive kernel and
satisfies $k^{\prime}\in L^1(R^+)$, then for  any $y\in L^1(R^+)$,
$$
\int_{t_0}^{t_1}|k*y(\tau)|^2d\tau 
\le   \beta_0K\int_{t_0}^{t_1}k*y \; \; y(\tau)d\tau, 
$$
where $K = |k|_1^2 + 4|k^{\prime}|_1^2$, $0 \le t_0 <
t_1 < \infty$, and $\beta_0 > 0$
is such that   $k(t) - \beta_0e^{-t} > 0$. Moreover, 
if $k(t) - \beta_0e^{-\gamma t} > 0 $ for some constant   
$\gamma >0$, then 
the conclusion also holds with $\beta_0 K$ replaced by
$\frac{1}{\gamma}\beta_0 K$.     
}

\bigskip

Now we estimate the norm of the concentration $C(x,t)$.

Taking the inner product  of (\ref{eqn}) with $C$ in $L^2(R)$ 
and integrating by
patrs, we obtain

\begin{eqnarray} 
\frac12 \frac{d}{dt}\|C\|^2 + D\|C_x\|^2 - 
\frac{1}{2}\int_{-\infty}^{\infty}
U^{\prime}(x)C^2 
dx + k(0)\|C\|^2 =   \nonumber  \\
 - \int_{-\infty}^{\infty}\int_0^t k^{\prime}(t-s)C(x, s)ds C(x,t)dx.
 \label{estimate1}
\end{eqnarray}
We assume that 
\begin{eqnarray}
\frac12 U^{\prime}(x) \le - \alpha_0 + k(0), \; \mbox{for some positive
                                        constant $\alpha_0$},
        \label{condition1}       \\
 k^{\prime} > 0, \mbox{and} \;\; k, k', k'' \in L^1(R^+),
        \label{condition2}
\end{eqnarray}
which imply  that the velocity field $U(x)$ has bounded gradient
and the probability distribution
of the random stopping times  is 
strictly increasing.

Integrating (\ref{estimate1}) with respect to $t$ 
from $t_0$ to $t_1$, and using (\ref{condition1}), (\ref{condition2})
and the Lemma above, we obtain that
\begin{eqnarray} 
\frac{1}{2}(\|C(t_1)\|^2 - \|C(t_0)\|^2) + \alpha_0\int_{t_0}^{t_1} 
\|C(s)\|^2 ds \le 0. 
\label{estimate2}
\end{eqnarray} 
Thus we have
$$\|C(t_1)\|^2 \le \|C(t_0)\|^2, \;\;\forall\;\;\;t_0 \le t_1.$$
If we take $t_0 = t, t_1 = t + 1$ in (\ref{estimate2}), then
\begin{eqnarray} 
\|C(t + 1)\|^2 
\le    \int_{t}^{t + 1} \|C(s)\|^2 ds 
\le    \frac{1}{2\alpha_0} (\|C(t)\|^2 - \|C(t + 1)\|^2),
\end{eqnarray} 
and
\begin{eqnarray*} 
\max\limits_{t\le s\le t+1}\|C(s)\|^2 
= \|C(t)\|^2 = \|C(t)\|^2 - \|C(t+1)\|^2 + \|C(t+1)\|^2   \\
\le (1 + \frac{1}{2\alpha_0}) 
(\|C(t)\|^2 - \|C(t + 1)\|^2).
\end{eqnarray*} 
By the difference inequality of Nakao \cite{nakao1} \cite{nakao2}, 
we finally obtain
\begin{eqnarray} 
\|C(t)\|^2 \le A e^{- B t}, \;t \ge 0, 
\end{eqnarray} 
for some positive constants $A$ and $B$. Thus $\|C(t)\|$
decays exponentially fast.
Therefore we have the following result.

\begin{theorem}
If the velocity field $U(x)$ has bounded gradient and 
the probability distribution $k(t)$ 
  (of the passive tracer random stopping times)
    is strictly increasing,
  namely, they satisfy conditions (\ref{condition1}) 
  and (\ref{condition2}), respectively,
  then the concentration of the passive tracer 
  with random stopping  approaches zero exponentially fast
  in mean-square norm.
\end{theorem}

\section{Exponential Decay  with History Source Term}

If there is a source term due to the past history of 
$C$ from $- \infty$ to $0$,
more precisely, if there is a source term  
$$ f(x, t) =  - \int_{-\infty}^0k'(t-s)C(x, s) dx, \;x \in R, t \ge 0,$$
then, (\ref{eqn2}) can be written as 
\begin{eqnarray}
C_t + U(x)C_x = DC_{xx} - k(0)C - \int_{-\infty}^tk'(t-s)C(x, s)ds
                        \nonumber  \\
 = DC_{xx} - k(0)C - \int_{0}^{\infty}k'(s)C(x, t - s)ds. 
 \label{eqn3}
\end{eqnarray}
Let 
\begin{eqnarray}
C^t(x, s) = C(x, t -s),  \\
\eta^t(x, s) = \int_0^s C^t(x, \tau)d\tau 
= \int_{t-s}^tC(x, \tau)d\tau, \; s \ge 0. 
\end{eqnarray}

We further assume that $\lim\limits_{s\to\infty}k'(s) = 0$, then

$$
 \int_{0}^{\infty}k'(s)C(x, t - s)ds = \int_{0}^{\infty}k'(s)C^t(x, s)ds 
= -
\int_{0}^{\infty}k''(s)\eta_t(s)ds.
$$
 Let $\mu(s) = k''(s)$, then (\ref{eqn3}) can be reformulated as
\begin{eqnarray} 
C_t + U(x)C_x = DC_{xx} - k(0)C - \int_{0}^{\infty}\mu(s)\eta^t(x,s)ds,  
        \label{eqn4}      \\
\eta_t^t(x, s) = C(x, t) - \frac{\partial}{\partial s}\eta^t. 
        \label{eqn5}
\end{eqnarray}
The initial   conditions are then
\begin{eqnarray} 
C(x, 0) = C_0(x), \;\;\eta^0(x, s) = \eta_0(x, s) 
   = \int_{-s}^0C(x, \tau)d\tau, \;x\in R. 
\label{IC2}
\end{eqnarray} 
Let
$$H = L^2(R) \times L^2(R^+; \mu; L^2(R)),
$$
where 
$$
L^2(R^+; \mu; L^2(R)) = 
\{\phi:R^+ \to L^2(R)|\int_0^{\infty}\mu(s)\|\phi(s)\|^2ds < \infty\} ,
$$ 
and $<\cdot>_{\mu}$, $\|\cdot\|_{\mu}$  denote  the inner product and 
norm in this space, respectively.

The well-posedness of the initial boundary value problem of 
(\ref{eqn4}), (\ref{eqn5})  and (\ref{IC2}) can be  shown as in
\cite{gr} and \cite{mi}. Here, we focus on the asymptotic behavior under 
the
condition  (\ref{condition1}) on the velocity field $U(x)$,
and the following conditions on the probability distribution $k(t)$
of the random stopping time 
\begin{eqnarray} 
k'> 0  \; \mbox{and} \; \; \lim\limits_{s\to\infty}k'(s) = 0, 
        \label{condition3}   \\
k'' \in C^1(R^+)\cap L^1(R^+)\;\; \mbox{and}\;\; 
        k''(s) \le 0, k'''(s) \ge 0,  \forall s \in R^+,  
        \label{condition4}              \\
k'''(s) + \delta k''(s)  \le  0, 
        \;\;\forall s \in R^+, \;\mbox{and for some constant}\; \delta > 0. 
                        \label{condition5}
\end{eqnarray} 

We now estimate the evolution of $C, \eta^t$ described by 
(\ref{eqn4})-(\ref{eqn5}).
Multiplying  (\ref{eqn4}) by $C$ and  (\ref{eqn5}) by $\eta^t$, then 
integrating by parts and adding the results, we get
$$
 \frac12 \frac{d}{dt}(\|C\|^2 + \|\eta^t\|_{\mu}^2) 
 - \int_{-\infty}^{\infty}U'(x)|C|^2dx $$ 
 \begin{eqnarray}
 + k(0)\|C\|^2 + D\|C_x\|^2
 + <\frac{\partial}{\partial s}\eta^t, \eta^t>_{\mu} = 0.  
 \label{eqn6}
\end{eqnarray}
Since for a fixed $t\in(0, \infty)$ both $\mu\|\eta^t\|$ and
$\mu\|\frac{\partial}{\partial s}\eta^t\|$ belong to $L^1(R^+)$,   
we now have
\begin{eqnarray}
 \lim\limits_{s\to\infty}\mu(s) \| \eta^t(s) \|^2 = 0. 
 \label{eqn7}
\end{eqnarray} 
Hence, using (\ref{eqn7}) and integrating by parts, we obtain
\begin{eqnarray} 
<\frac{\partial}{\partial s}\eta^t, \eta^t>_{\mu} = -
\int_0^{\infty}\mu'(s)\|\eta^t\|^2 ds. 
\label{eqn8}
\end{eqnarray}
By putting (\ref{condition1}), (\ref{condition3})-(\ref{condition5}), 
and
(\ref{eqn8}) into  (\ref{eqn6}),   we finally get
\begin{eqnarray}
 \frac{d}{dt}(\|C\|^2 + \|\eta^t\|_{\mu}^2) + 2\alpha_0\|C\|^2 +
        2\delta\|\eta^t\|_{\mu}^2 \le 0,
\end{eqnarray}
where $\alpha_0, \delta$ are defined in (\ref{condition1}), 
(\ref{condition5}),
respectively.
So by the Gronwall inequality \cite{Temam}, 
we obtain the exponential decay for $\|C\|^2 + \|\eta^t\|_{\mu}^2$
and thus for the concentration (in mean-square norm) $\|C\| $. 
We summarize this result in the following theorem.

\begin{theorem}
 Assume that the  velocity field $U(x)$ has bounded gradient and the
  probability distribution $k(t)$ 
  (of the passive tracer random stopping times) is   strictly 
increasing,
  namely, they satisfy conditions (\ref{condition1}) 
  and (\ref{condition2}), respectively. If $k(t)$ satisfies
  further growth conditions (\ref{condition3}), (\ref{condition4})
  and (\ref{condition5}), 
  then the concentration of the passive tracer 
  with random stopping and with a history source 
   approaches zero exponentially fast
  in mean-square norm.  
\end{theorem}

\vspace{0.2cm}
Note that the transformation that leads to
the reformulation (\ref{eqn4})-(\ref{eqn5}) is inspired by \cite{gpm}. 
Under this transformation, a  nonautonomous system (\ref{eqn3})
   with  history source term   becomes an autonomous sysytem 
(\ref{eqn4})-(\ref{eqn5}), at the  expense of increasing 
system components or dimension.

\section{Discussions}

In this paper, we have studied the nonlocal behavior of 
passive tracer dispersion with random stopping, using an
advection-diffusion equation with a memory or history term.
 We have shown the exponential decay of the
passive tracer concentration, under suitable conditions for the
velocity field and the probability distribution of random stopping time.
There are other situations where tracer  dispersion
is modeled by advection-diffusion type of equations
with memory terms \cite{Smith}, \cite{Davis}
\cite{Koch} and \cite{Cushman}.

\bigskip

{\bf Acknowledgement.} This work was supported by the U. S. National 
Science
                Foundation Grant DMS-9704345, the National  
                Natural Science  Foundation  of China Grant   19701023,
                and  by the Clemson Center for Industrial Mathematics 
and
                Statistics.

\end{document}